\documentclass[11pt]{article}

\usepackage{geometry}
\usepackage{graphicx}
\usepackage{amsmath}
\usepackage{color}
\usepackage{diagbox}
\usepackage{amssymb}
\usepackage{booktabs}
\usepackage{array}
\usepackage{paralist}
\usepackage{verbatim}
\usepackage{subfig}
\usepackage{indentfirst}
\usepackage{fancyhdr}
\title{On Multiplayer Chomp Theory}
\author{Purui Zhang\footnote{School of Mathematical Sciences, Fudan University, 220 Handan Road, Yangpu District, Shanghai, China. 200433\quad boshi\_haha@126.com}}
\date{}

\begin{document}
\maketitle
\paragraph{Abstract} Decades after David Gale\cite{1} presented the concept of Chomp and S.-Y.R. Li\cite{2} produced his very first multiplayer model to investigate Multiplayer Nim, we hereby establish another Multiplayer Model to specifically analyze Chomp. Under such model, we will re-define `positions', `rules' and other basic concepts and try our best to analyze the relationship between different rules and corresponding positions.
\tableofcontents\clearpage
\section{Chomp Positions}
\subsection{Positions and Chomp Move}
In standard two-player chomp game, what we care about is the winner (and the loser) of a certain position. To be specific, given a starting position (for example, an $m$ by $n$ chocolate bar), player 1 and 2 take turns to `chomp' it, and the player who take the last bite loses. Because no matter how many players there would be, the `chomp' move for all players are same. So it is intuitive to define positions and moves firstly for better setting up multiplayer rules.
\\~\\
\paragraph{Definition 1.1} A {\bfseries position} is a tuple $(a_1,a_2,...,a_n)$ with $a_1\geqslant a_2\geqslant...\geqslant a_n\geqslant0$. All elements are non-negative integers. Also, positions satisfy that adding or deleting `0's at the end remain the position unchanged. For example, $(3,2,0)$ and $(3,2)$ are same positions. Also we define the volume of position $P$ as $|P|=\sum\limits_{i=1}^{n}a_i$.\\~\\
For $m$ by $n$ rectangular chocolate bar, $(a_1,a_2,a_3,...,a_n)$ with $a_1=a_2=...=a_n=m$ and $(a_1,a_2,a_3,...,a_m)$ with $a_1=a_2=...=a_m=n$ both represent it.
\begin{table}[!htbp]
\centering
\begin{tabular}{lllll}
x & x & x & x & x \\
x & x & x & x & x \\
x & x & x & x & x \\
x & x & x & x & x
\end{tabular}
\caption{Position $(4,4,4,4,4)$ representing a 4 by 5 chocolate.}
\label{rect}
\end{table}
Empty position $(0)$ represents the ending state of chomp, with no chocolate pieces left.\\~\\
The volume of position can be interpreted as the number of chocolate pieces in the given chocolate bar. \\~\\
\paragraph{Definition 1.2} A chomp move is an operation on position $P=(a_1,a_2,...,a_n)\to Q=(b_1,b_2,...,b_n)$ such that:
\begin{itemize}
    \item given $1\leqslant x\leqslant n$ and non-negative integer $a$;
    \item for any $1\leqslant i\leqslant x$, $a_i=b_i$;
    \item for any $x\leqslant j\leqslant n$, $b_j=a$ if $a_j>a$ and $b_j=a_j$ if $a_j\leqslant a$;
    \item $P$ and $Q$ are different positions.
\end{itemize}

For example, we know that $(5,3,3)\to(5,3,2)\to(1,1,1)\to(1,0,0)\to(0)$ are valid, consecutive chomp moves. And also obviously, for volume, we have

\paragraph{Fact 1.3} If $P\to Q$, then $|P|>|Q|$.

\subsection{Position Sets}

\paragraph{Definition 1.4} Let position $P=(a_1,a_2,...,a_n)$ where $a_1\geqslant a_2\geqslant...\geqslant a_n\geqslant0$. Define
\begin{itemize}
    \item $\mathcal{L}=\{P|P\neq(0)\}$;
    \item $\mathcal{L}_0=\{P|P=(1)\}$;
    \item $\mathcal{L}_1=\{P|P=(1,1,1...,1)\text{ or }(a_1)\}\backslash\mathcal{L}_0$;
    \item $\mathcal{L}_2=\mathcal{L}\backslash(\mathcal{L}_0\bigcup\mathcal{L}_1)$;
    \item $\mathcal{D}_k=\{P|\min\{a_1,n\}=k\}$
\end{itemize}
We could easily get
\begin{itemize}
    \item $\mathcal{L}=\biguplus\limits_{i=0}^{2}\mathcal{L}_i=\biguplus\limits_{k=1}^{\infty}\mathcal{D}_k$;
    \item $\mathcal{D}_1=\mathcal{L}_0\bigcup\mathcal{L}_1$
\end{itemize}
Specifically, to understand why we define $\mathcal{L}_0$, $\mathcal{L}_1$, and $\mathcal{L}_2$, we investigate the `shape' of positions. Intuitively, a 1 by 1 chocolate piece is like a 0-dimensional `point', and 1 by n ($n>1$) chocolate bar is like a 1-dimensional `line'.\\~\\
\paragraph{Fact 1.5} If $P\to Q$,
\begin{itemize}
    \item $Q\in\mathcal{L}_0\Rightarrow P\notin\mathcal{L}_2$;
    \item $P\in\mathcal{D}_p$, $Q\in\mathcal{D}_q\Rightarrow p\geqslant q$
\end{itemize}
In many cases we investigate positions that can be chomped from, or can chomp to a given position. So we can further define the set of possible moves as follows.

\paragraph{Definition 1.6} Given position $P$, define
\begin{itemize}
    \item $\text{Mov}P=\{Q|P\to Q\}$;
    \item $\overline{\text{Mov}}P=\{Q|Q\to P\}$
\end{itemize}
\paragraph{Fact 1.7}
\begin{itemize}
    \item $P\to Q\Leftrightarrow Q\in\text{Mov}P\Leftrightarrow P\in\overline{\text{Mov}}Q$;
    \item $\overline{\text{Mov}}(0)=\mathcal{L}$;
    \item $|\text{Mov}P|=|P|$ because in the chocolate bar corresponding to position $P$, you can choose exactly $|P|$ different position to chomp with.
\end{itemize}
\section{Multiplayer Rule of Chomp}
\subsection{Player's Preferences}
After we set up the positions and possible moves, we can now invite our players into the game. The easiest way to conduct the game is to let player 1 through $n$ chomp the position by turn, and then we evaluate their score according to the player who take the last move. \\~\\
In 2-player game, the one who takes the last bite loses. We can also say that the one who takes the last bite gets 0 points, and another one gets 1 point. To generalize, we predetermine a `score sequence' $(\alpha_1\alpha_2...\alpha_n)$ for $n$ players, and we make the rule that the one who takes the last bite (assume that is player $i$) gets $\alpha_1$ points, and the one before (player $i-1$) it gets $\alpha_2$ points, and the one before that (player $i-2$) gets $\alpha_3$ points, etc. According to this process, we can set up the concept of multiplayer rules. \\~\\
\paragraph{Definition 2.1} A multiplayer rule of chomp is represented as $f=(\alpha_1\alpha_2...\alpha_n)$, where $\alpha_1$ through $\alpha_n$ are distinct real numbers and $n$ indicates the number of players (also denoted as $|f|$). Without special indication, we assume $n\geqslant 2$. When playing multiplayer chomp under such rule,
\begin{itemize}
    \item the player who takes the last bite (assume that is player $i$) gets $\alpha_1$ points, and the one before (player $i-1$) it gets $\alpha_2$ points, and the one before that (player $i-2$) gets $\alpha_3$ points, etc.
    \item All players choose moves that can help them achieve the highest point. If there are multiple solutions, pick one of them randomly;
    \item If chomp to (0) is one of the most beneficial moves, the player must do this to avoid long rally (this entry is also called `player's preference').
\end{itemize}

Player's preference does not show up in the well-known 2-player chomp, but it is necessary to simplify our further analysis of multiplayer rules.\\ \\

\subsection{Zermelo's Theorem of Multiplayer Chomp}

When rules are given, multiplayer chomp is a combinatorial game with perfect information for all players. We believe that two-player Zermelo's Theorem can be extended to multiplayer chomp. \\~\\
\paragraph{Theorem 2.2} Given $n$-player rule $(\alpha_1\alpha_2...\alpha_n)$. Then for any $P\in\mathcal{L}$ as starting position, the score player 1 could get is determined.\\~\\
\textit{proof.}\quad We use induction on $|P|$. \\~\\
If $|P|=1$, then $P=(1)$. Player 1 can only chomp it to (0), which means that player 1 must get $\alpha_1$ points;
If for any position with $|P|\leqslant k=1$ the score of player 1 has been determined, we consider position with $|P|=k$. \\~\\
For any $Q\in\text{Mov}P$, $|Q|<|P|=k$ holds, so the score of player 1 when $Q$ is the starting position has been determined. Let $|\text{Mov}P|=t\leqslant k-1$ and $\text{Mov}P={Q_0, Q_1, ..., Q_{t-1}}$ (with $Q_0=(0)$). We further assume that the score player 1 could get when $Q_i$ is the starting position is $\alpha_{q_i}$ ($0\leqslant i\leqslant t-1$). Then, if player 1 on starting position $P$ does the chomp move $P\to Q_i$, player 2 would get exactly $\alpha_{q_i}$ points. That means that player 1 would get $\alpha_{q_i+1}$ points (if $q_i=n$ then let $q_i+1=1$). According to player 1's preferences, $\max\limits_{0\leqslant i\leqslant t-1}\{\alpha_{q_i+1}\}$ would be chosen for player 1's highest score. Because the set above is finite, maximum can always be reached. Thus, according to induction, theorem holds. $\square$
\subsection{Rules and Solutions}
With respect to the proof above, we can see that if player 1 plays different chomp moves at the beginning, the outcome of player 1's score varies, but only the most beneficial (also only choose (0) if its the most beneficial) moves will be chosen. \\~\\
\paragraph{Definition 2.3} For rule $f$, any $P\in\mathcal{L}$, define the solution set $$\text{Sol}_fP=\{Q|P\to Q\text{ and player 1 scores highest with this move, obeying player's preference}\}$$ and
$$\overline{\text{Sol}}_fP=\{Q|Q\to P\text{ and player 1 scores highest with this move, obeying player's preference}\}$$\\~\\
What's more, if $Q\in\text{Sol}_fP$, we call that $Q$ is a solution of $P$, or $P\to Q$ is a solution move. \\~\\
\paragraph{Fact 2.4}For any rule $f$, apparently we have $\text{Sol}_fP\subset\text{Mov}P$ and $\overline{\text{Sol}}_fP\subset\overline{\text{Mov}}P$ for any $P\in\mathcal{L}$.\\~\\
\paragraph{Definition 2.5} For rule $f$, call $P_1\to P_2\to P_3\to...\to P_k\to P_{k+1}=(0)$ a solution chain if for any $1\leqslant i\leqslant k$, $P_i\to P_{i+1}$ is a solution move under rule $f$.
\subsection{Ordinals}
Indeed under certain rules, the set of solutions and the highest score player 1 would get completely reflects the properties of a position. However, these are not the most essential ones. The concept which we will introduce later has a pivotal role in multiplayer chomp rule analysis, and all our further research cannot proceed without it.\\~\\
\paragraph{Definition 2.6}For $n$-player game with rule $f=(\alpha_1\alpha_2...\alpha_n)$ and given position $P\in\mathcal{L}$, if player 1 gets the highest score $\alpha_i$ when $P$ is the starting position, we define the ordinal of position $P$ as $\text{ord}_fP=i$. Also define $\text{ord}_f(0)=0$.\\~\\
So we know that all positions including (0) have their ordinals under certain rules. and all ordinals have the maximum of $n$ and minimum of 0. Now we can investigate more properties of ordinals. \\~\\
\paragraph{Theorem 2.7}If $P\to Q$ is a solution move under rule $f$, then $\text{ord}_fQ\leqslant n-1$ and $\text{ord}_fP=\text{ord}_fQ+1$.\\~\\
\textit{proof.}\quad Firstly, if $\text{ord}_fQ=n$, then $Q\neq(0)$ because $n$ (as the total number of players) must be a positive integer. If $Q$ is the starting position, player 1 would get a point of $\alpha_{\text{ord}_fQ}=\alpha_n$. So, if $P$ is the starting position, and $P\to Q$ is a solution move, player 2 would get the point $\alpha_{n}$, which makes the conclusion that player 1 would get $\alpha_{1}$ points. However, $P\to(0)$ also make player 1 get exactly $\alpha_1$ points, so according to player's preferences, $P\to(0)$ is the only solution and $P\to Q$ is not. This leads to a contradiction. Thus, $\text{ord}_fQ\leqslant n-1$.\\~\\
For the latter half of the theorem, also assume that if $Q$ is the starting position, player 1 would get a point of $\alpha_{\text{ord}_fQ}$. So, if $P$ is the starting position, and $P\to Q$ is a solution move, player 2 would get the point $\alpha_{\text{ord}_fQ}$, which makes the conclusion that player 1 now would get $\alpha_{\text{ord}_fQ+1}$ points. Given that $\text{ord}_fQ+1\leqslant n$, the ordinal of $P$ is $\text{ord}_fP=\text{ord}_fQ+1$. $\square$\\~\\
\paragraph{Theorem 2.8} For $n$-player game with rule $f=(\alpha_1\alpha_2...\alpha_n)$, if we have solution chain $P_1\to P_2\to P_3\to...\to P_k\to P_{k+1}=(0)$, then $k\leqslant n$.\\~\\
\textit{proof.}\quad Using theorem 2.7 $k$ times we can get $\text{ord}_fP_i=k+1-i$ for $1\leqslant i\leqslant k$. When $i=1$ we have $\text{ord}_fP_1=k$, so obviously $k\leqslant n$. $\square$\\~\\
According to theorem 2.7 and 2.8 we can get more quick facts. \\~\\
\paragraph{Fact 2.9}For $n$-player rule $f$ and position $P$,
\begin{itemize}
    \item $\text{ord}_fP=n\Leftrightarrow\overline{\text{Sol}}P=\varnothing$;
    \item $\text{ord}_fP=i\enspace(i<n)\Leftrightarrow\forall Q\in\overline{\text{Sol}}P,\enspace\text{ord}_fQ=i+1$;
    \item $\text{ord}_fP=1\Leftrightarrow\text{Sol}P=\{(0)\}$;
    \item $\text{ord}_fP=i\enspace(i>1)\Leftrightarrow\forall Q\in\text{Sol}P,\enspace\text{ord}_fQ=i-1$
\end{itemize}
After investigating the properties of ordinals, we can construct a way to calculate ordinals of positions using recurrence relation: \\~\\
Consider $n$-player rule $f=(\alpha_1\alpha_2...\alpha_n)$ and starting position $P\in\mathcal{L}$, to find out the score player 1 could get (that is, $\alpha_{\text{ord}_fP}$), we first list all possible moves from $P$. Then, for any $Q\in\text{Mov}P$, we take that as the starting position and see what player 1 could get from that (which is in fact $\alpha_{\text{ord}_fQ}$). Thus, as we go back analyze position $P$ as starting position, we know that if player 1 does $P\to Q$, player 2 would also get exactly $\alpha_{\text{ord}_fQ}$. If $\text{ord}_fQ=n$, $P\to Q$ is never a solution move due to the first fact of fact 2.9. If we rule out those positions and only leave with  $\text{ord}_fQ\leqslant n-1$, we yield that player 1 would get $\alpha_{\text{ord}_fQ+1}$ points. Taking the maximum score of all $Q\in\text{Mov}P$ satisfying $\text{ord}_fQ\neq n$ leads us to the highest score for player 1 with starting position $P$. Hence we get the fact that\\~\\
\paragraph{Fact 2.10} For $n$-player rule $f=(\alpha_1\alpha_2...\alpha_n)$ and position $P\in\mathcal{L}$,
$$\alpha_{\text{ord}_fP}=\max\limits_{Q\in\text{Mov}P,\enspace\text{ord}_fQ\neq n}\{\alpha_{\text{ord}_fQ+1}\}$$
\paragraph{Definition 2.11} Given rule $f$, positive integer $i$ and position $P$, define resolvent set $$\text{Res}_f^iP=\{Q|P\to Q, \text{ord}_fQ=i-1\}$$
With this supporting definition, we immediately get $\text{ord}_fP=i\Leftrightarrow\text{Res}_f^iP=\text{Sol}_fP$. Fact 2.10 can be simplified as $$\alpha_{\text{ord}_fP}=\max\limits_{i=1}^{|f|-1}\{\alpha_i|\text{Res}_f^iP\neq\varnothing\}$$
\paragraph{Corollary 2.12} If $P=(1)$, then for any rule $f$, $\text{Res}_f^1P=\{(0)\}$; for any positive integer $i\geqslant2$, $\text{Res}_f^iP=\varnothing$.
\section{Rule Isomorphism}
In standard 2-player chomp, our ultimate goal is to judge whether the first player wins or loses on any given position. We could have started the same investigation on each of the multiplayer rules; however, not all rules are needed to be analyzed. In fact, similarities can be found in many different rules.
\subsection{Definition of Isomorphism}
Given a rule and a position, we can decide all solutions, the score, and the ordinal. The solution set and the ordinals are very important because they help us to sketch out the overall structure that the rule forms. Hence, it is intuitive to carry out the concept of rule isomorphism. \\~\\
\paragraph{Theorem 3.1}For two rules $f$ and $g$ (do not need to have same number of players),
$$\forall P\in\mathcal{L}, \text{Sol}_fP=\text{Sol}_gP\Leftrightarrow\forall P\in\mathcal{L}, \text{ord}_fP=\text{ord}_gP$$
\textit{proof.}\quad First prove $\text{Sol}_fP=\text{Sol}_gP\Rightarrow\text{ord}_fP=\text{ord}_gP$. We can use induction on $|P|$. If $|P|=1$, then $P=(1)$, $\text{Sol}_fP=\text{Sol}_gP=\{(0)\}$. So $\text{ord}_fP=\text{ord}_gP=0+1=1$. If for any $|P|<k$ where $k\geqslant2$, theorem holds, we investigate a position $P$ with $|P|=k$. Because $\text{Sol}_fP=\text{Sol}_gP$, we can pick position $Q$ so that $P\to Q$ is the solution move of both rule $f$ and $g$. According to fact 1.3 and theorem 2.7, $|Q|<|P|=k$, yielding $\text{ord}_fP=\text{ord}_fQ+1=\text{ord}_gQ+1=\text{ord}_gP$. According to induction, $\text{ord}_fP=\text{ord}_gP$ holds for all $P\in\mathcal{L}$.\\~\\
Then we prove the inverse proposition. We also use the induction on $|P|$. If $|P|=1$, then $P=(1)$, $\text{Sol}_fP=\text{Sol}_gP=\{(0)\}$. If for any $|P|<k$ where $k\geqslant2$, theorem holds, we investigate a position $P$ with $|P|=k$. If $\text{ord}_fP=\text{ord}_gP=1$, then according to player's preference, we have $\text{Sol}_fP=\text{Sol}_gP=\{(0)\}$; if $\text{ord}_fP=\text{ord}_gP>1$, we know that $$\text{Sol}_fP=\{Q|Q\in\text{Mov}P, \text{ord}_fQ=i-1\}$$
$$\text{Sol}_gP=\{Q|Q\in\text{Mov}P, \text{ord}_gQ=i-1\}$$
So obviously $\text{Sol}_fP=\text{Sol}_gP$. According to induction, $\text{Sol}_fP=\text{Sol}_gP$ holds for all $P\in\mathcal{L}$.$\square$\\~\\
\paragraph{Definition 3.2}\quad If two equations in theorem 3.1 holds, denote $f\sim g$, calling two rules isomorphs. $\sim$ is apparently an equivalence relation.\\~\\
Recall that for a rule $f=(\alpha_1\alpha_2\alpha_3...\alpha_n)$, we only need $\alpha_1$ through $\alpha_n$ to be distinct real numbers. Hence, we need to restrict the range of those numbers so that we can better compare them. \\~\\
{\bfseries Theorem 3.3}\quad Given rule $f=(\alpha_1\alpha_2\alpha_3...\alpha_n)$ and $g=(\beta_1\beta_2\beta_3...\beta_n)$, if for any $i\neq j$, $(\alpha_i-\alpha_j)(\beta_i-\beta_j)>0$, then $f\sim g$.\\~\\
\textit{proof.}\quad We need to prove that for any position $P\in\mathcal{L}$, $\text{ord}_fP=\text{ord}_gP$. We use induction on $|P|$. \\~\\
If $|P|=1$, then $P=(1)$. $\text{ord}_fP=\text{ord}_gP=1$. If for any $|P|\leqslant t-1$, $\text{ord}_fP=\text{ord}_gP$, then we look at position $P$ with $|P|=t$. According to Fact 2.10, $$\alpha_{\text{ord}_fP}=\max\limits_{Q\in\text{Mov}P,\enspace\text{ord}_fQ\neq n}\{\alpha_{\text{ord}_fQ+1}\}$$
$$\beta_{\text{ord}_gP}=\max\limits_{Q\in\text{Mov}P,\enspace\text{ord}_gQ\neq n}\{\beta_{\text{ord}_gQ+1}\}$$
Let $$\mathcal{P}_f=\{Q|Q\in\text{Mov}P,\enspace\text{ord}_fQ\neq n\}$$
$$\mathcal{P}_g=\{Q|Q\in\text{Mov}P,\enspace\text{ord}_gQ\neq n\}$$
With respect to our induction hypothesis, $\mathcal{P}_f=\mathcal{P}_g:=\mathcal{P}$. For any position $R\in\mathcal{P}$ we also know that $\text{ord}_fR=\text{ord}_gR$.\\~\\
Now we can finally assert that $\text{ord}_fP=\text{ord}_gP$. If not, because $\alpha_{\text{ord}_fP}=\max\limits_{Q\in\mathcal{P}}\{\alpha_{\text{ord}_fQ+1}\}$ and $\beta_{\text{ord}_gP}=\max\limits_{Q\in\mathcal{P}}\{\beta_{\text{ord}_gQ+1}\}$, we can find the maxima such that $\alpha_{\text{ord}_fP}=\alpha_{\text{ord}_fQ_1+1}$ and $\beta_{\text{ord}_gP}=\beta_{\text{ord}_gQ_2+1}\}$. Under rule $f$, we have $\alpha_{\text{ord}_fQ_1+1}>\alpha_{\text{ord}_fQ_2+1}$; but under rule $g$, we have $\beta_{\text{ord}_gQ_1+1}<\beta_{\text{ord}_gQ_2+1}$. Given that $\text{ord}_fQ_1+1=\text{ord}_gQ_1+1$, $\text{ord}_fQ_2+1=\text{ord}_gQ_2+1$, we finally have $(\alpha_{\text{ord}_fQ_1+1}-\alpha_{\text{ord}_fQ_2+1})(\beta_{\text{ord}_fQ_1+1}-\beta_{\text{ord}_fQ_2+1})<0$, leading to contradiction.\\~\\
According to induction, we successfully prove the theorem. $\square$\\~\\
With theorem 3.3, we can `normalize' all positions: previously, rule $f=(\alpha_1\alpha_2...\alpha_n)$ is with real number $\alpha_1,\enspace\alpha_2\enspace...\alpha_n$; but this time we can add an extra restriction that $(\alpha_1\alpha_2...\alpha_n)$ is a permutation of $(0123...n-1)$. Apparently,\\~\\
{\bfseries Fact 3.4}\quad For any rule $f=(\alpha_1\alpha_2...\alpha_n)$ there is only one `normalized' rule $(\beta_1\beta_2...\beta_n)$ (which is a permutation of $(0123...n-1)$) such that for any $i\neq j$, $(\alpha_i-\alpha_j)(\beta_i-\beta_j)>0$. $(\beta_1\beta_2...\beta_n)$ is called the normalized rule of $f$.\\ \\
From now on, all rules are presumed normalized unless specified.
\subsection{Simple Rules}
There are some kind of rules, say (687351024), only has small ordinals (in this example, you will find that the maximum ordinal for all positions is only 3). That means, game will go to an end within only a few moves. Given the tool of rule isomorphism, can we reduce it to a simpler rule and investigate that instead? The answer is yes. \\~\\
{\bfseries Theorem 3.5}\quad (Simple Rule Theorem) For $n$-player rule $f=(\alpha_1\alpha_2...\alpha_n)$, if for any $P\in\mathcal{L}$, $\text{ord}_fP\leqslant m<n$, let $g=(\alpha_1\alpha_2...\alpha_m)$, then $f\sim g$.\\~\\
\textit{proof.}\quad We need to prove that for any position $P\in\mathcal{L}$, $\text{ord}_fP=\text{ord}_gP$. We use induction on $|P|$. \\~\\
If $|P|=1$, then $P=(1)$. $\text{ord}_fP=\text{ord}_gP=1$. If for any $|P|\leqslant t-1$, $\text{ord}_fP=\text{ord}_gP$, then we look at position $P$ with $|P|=t$. Apparently, if $P\to Q$ is a solution move under rule $f$, we have $\text{ord}_fQ\leqslant m-1$ because $\text{ord}_fP\leqslant m$. Then, according to Fact 2.10,
$$\alpha_{\text{ord}_fP}=\max\limits_{Q\in\text{Mov}P,\enspace\text{ord}_fQ\leqslant m-1}\{\alpha_{\text{ord}_fQ+1}\}$$
$$\alpha_{\text{ord}_gP}=\max\limits_{Q\in\text{Mov}P,\enspace\text{ord}_gQ\leqslant m-1}\{\alpha_{\text{ord}_gQ+1}\}$$
Let $$\mathcal{P}=\{Q|Q\in\text{Mov}P,\enspace\text{ord}_fQ\leqslant m-1\}$$
With respect to our induction hypothesis, for any position $R\in\mathcal{P}$ we also know that $\text{ord}_fR=\text{ord}_gR$. Remember that we need to prove $\text{ord}_fP=\text{ord}_gP$. Now we assume that $\text{ord}_fP\neq\text{ord}_gP$ and we want to deduce a contradiction. Then according to Fact 2.10 there exists  maxima such that $\text{ord}_fP=\text{ord}_fQ_1+1$ and $\text{ord}_gP=\text{ord}_gQ_2+1$. Because $Q_1,\enspace Q_2\in\mathcal{P}$, $\alpha_{\text{ord}_fQ_1}>\alpha_{\text{ord}_fQ_2}$ and $\alpha_{\text{ord}_fQ_2}>\alpha_{\text{ord}_fQ_1}$ holds simultaneously. This leads to a contradiction. \\~\\
With respect to induction, we finally prove the whole theorem. $\square$\\~\\
The theorem above is also called `Reduction Theorem', it can helps us to `reduce' so many seemingly-complicated rules to easy ones. To be specific, \\~\\
{\bfseries Definition 3.6}\quad For $n$-player rule $f$, if there exists a position $P\in\mathcal{L}$ such that $\text{ord}_fP=n$, then we call $f$ a simple rule. Else, it is called a non-simple rule. According to Theorem 3.5, each non-simple rule has a simple rule isomorph.
\subsection{Transposition Properties}
\begin{table}[!htbp]
\centering
\begin{tabular}{lllllllllllll}
x &   &   &   &  &  &  &  &   &   &   &   &   \\
x &   &   &   &  &  &  &  & x &   &   &   &   \\
x & x &   &   &  &  &  &  & x &   &   &   &   \\
x & x &   &   &  &  &  &  & x & x & x &   &   \\
x & x & x & x &  &  &  &  & x & x & x & x & x
\end{tabular}
\caption{Position $P=(5,2,1,1)$ and its transposition $P^T=(4,2,2,1,1)$.}
\label{transposition}
\end{table}
Given a piece of chocolate, we can `flip' it along diagonal. To describe the `flip' algebraically, we introduce the definition of position transposition. \\~\\
{\bfseries Definition 3.7}\quad Given position $P=(a_1,a_2,...,a_n)$, define its transposition $P^T=(b_1,b_2,...,b_m)$ by letting $m=a_1$ and $b_i=|\{j\in\{1,2,...,n\}|a_j\geqslant i\}|$ for any $1\leqslant i\leqslant m$. Apparently transposition has straightforward properties like $(P^T)^T=P$, $\text{Mov}P^T=\{Q^T|Q\in \text{Mov}P\}$, $P\to Q\Leftrightarrow P^T\to Q^T$ and $|P^T|=|P|$.\\~\\
It is intuitive to argue that transposition does not alter the core properties of a position under a given rule. And we know that under certain rule, ordinal is the best way to characterize a position, so we can further assert that \\~\\
{\bfseries Fact 3.8}\quad For any position $P\in\mathcal{L}$ and any rule $f=(\alpha_1\alpha_2...\alpha_n)$,
\begin{enumerate}
\item $\text{ord}_fP^T=\text{ord}_fP$;
\item $\text{Sol}_fP^T=\{Q^T|Q\in\text{Sol}_fP\}$.
\end{enumerate}

\subsection{Restricted Isomorphism}
Given rule $f$ and $g$, If $f\sim g$ is true, you may try to prove this according to the definition of rule isomorphism. However, what if we want to prove that $f\not\sim g$? A straightforward way is to find a position as an counterexample such that ordinals of $f$ and $g$ at that spot are different. Then, where could we find such positions? We need restricted isomorphism to make our search range smaller and to discover more patterns. \\~\\
{\bfseries Definition 3.8}\quad If for any position $P$ within position set $\mathcal{S}\subset\mathcal{L}$, $\text{ord}_fP=\text{ord}_gP$, then we call $f$ and $g$ are isomorphs on set $\mathcal{S}$, which is denoted as $f\overset{\mathcal{S}}{\sim}g$. Obviously, $f\sim g\Rightarrow f\overset{\mathcal{S}}{\sim}g$ and $f\overset{\mathcal{S}}{\not\sim}g\Rightarrow f\not\sim g$.\\~\\
Restricted isomorphism will be much easier to analyze because the `structure' of the restricted set is simpler than the whole $\mathcal{L}$. \\~\\
{\bfseries Question 3.9}\quad Explore the necessary and sufficient condition for $f\overset{\mathcal{D}_1}{\sim}g$.\\~\\
\textit{solution.}\quad There are two classes of positions in $\mathcal{D}_1$, one is $(1),(2),(3)...,(k),...$, another is $(1),(1,1),(1,1,1), ...$. We only need to look at the first class. It is easy to conclude that $(m)\to(n)$ if $m>n$. Given rule $f=(\alpha_1\alpha_2...\alpha_n)$, we can conclude that:
\begin{itemize}
\item $\text{ord}_f(1)=1$;
\item if $\alpha_1>\alpha_2$, $\text{ord}_f(2)=1$; else, $\text{ord}_f(2)=2$;
\item if $\alpha_1>\alpha_2$, $\text{ord}_f(3)=1$; else, if $\alpha_2>\alpha_3$, $\text{ord}_f(3)=2$; else, $\text{ord}_f(3)=3$;
\item ...
\end{itemize}
We hereby give some examples to better illustrate the process:
\begin{itemize}
\item $f=(01243)$, $0<1<2<4$, $4>3$, so $\text{ord}_f(1)=1$, $\text{ord}_f(2)=2$, $\text{ord}_f(3)=3$, $\text{ord}_f(4)=4$, $\text{ord}_f(5)=4$, $\text{ord}_f(6)=4$, ... note that if $(5)$ is starting position, $(5)\to(4)$ makes player 1 get 3 points, but $(5)\to(3)$ makes player 1 get 4 points. So $(5)\to(3)$ is solution move.
\item $f=(013245)$, $0<1<3$, $3>2$, so $\text{ord}_f(1)=1$, $\text{ord}_f(2)=2$, $\text{ord}_f(3)=3$, $\text{ord}_f(4)=3$, $\text{ord}_f(5)=3$, $\text{ord}_f(6)=3$, ...
\item $f=(0123)$, $0<1<2<3<4$, so $\text{ord}_f(1)=1$, $\text{ord}_f(2)=2$, $\text{ord}_f(3)=3$, $\text{ord}_f(4)=3$, $\text{ord}_f(5)=4$, $\text{ord}_f(6)=4$, ... note that if $(5)$ is starting position, $(5)\to(4)$ makes player 1 get 0 points, but $(5)\to(3)$ makes player 1 get 3 points. $(5)\to(3)$ is solution move.
\end{itemize}
So, we define $m_f=\min\{i\in\mathbb{Z}_+|\alpha_{i+1}-\alpha_i<0\}$ ($\alpha_{n+1}:=\alpha_1$ for convenience). It is now clear that
\begin{itemize}
\item $\text{ord}_f(i)=i$ if $i\leqslant m_f$;
\item $\text{ord}_f(i)=m_f$ if $i>m_f$.
\end{itemize}
Finally, we could get that $f\overset{\mathcal{D}_1}{\sim}g\Leftrightarrow m_f=m_g$. $\square$\\~\\
Also we have an important by-product: \\~\\
{\bfseries Corollary 3.10}\quad If given $i>m_f$, then $\overline{\text{Mov}}\{P|\text{ord}_fP=i\}\subset\mathcal{L}_2$.
\subsection{Standard Rules}
Now it is time to generalize our two-player chomp, with the rule representation (01).\\~\\
{\bfseries Definition 3.11}\quad Define $s_n=(0123...n-1)$ as an $n$-player simple standard rule. If another rule $g\sim s_n$, we call $g$ a standard rule. All non-standard rules do not have any simple standard rule isomorphs.\\~\\
According to the results given in question 3.9, we know that $$m_{s_n}=\min\{i\in\mathbb{Z}_+|\alpha_{i+1}-\alpha_i<0\}=n$$
Moreover, no other $n$-player rules satisfies such equation. So we can quickly get that $n$-player simple standard rules could never be isomorphs of any other rules with players no more than $n$. We hereby indicate an important property of non-standard rules.\\~\\
{\bfseries Theorem 3.12}\quad $f$ is a non-standard rule $\Leftrightarrow$ ($\text{ord}_fP=1\Leftrightarrow P=(1)$).\\~\\
\textit{proof.}\quad We only need to discuss simple rules as non-simples have simple rule isomorph. Assume $f=(\alpha_1\alpha_2\alpha_3...\alpha_n)$ is simple non-standard. For the only position $(1)\in\mathcal{L}_0$, $\text{ord}_fP=1$ is obvious; for any position $T\in\mathcal{L}_1$, as $m_f<n$, we know that $\text{ord}_fT<n$. Then, because for any position $R\in\mathcal{L}_2$, where $R=(a_1, a_2, a_3, ..., a_m)$, $m>1$, $a_1>2$, we know that $R\to Q=(a_1)\in\mathcal{L}_1$ is a legal move. So player 1 at position $R$ could get at least $\alpha_{\text{ord}_fQ+1}$ points. As $Q\in\mathcal{L}_1$ has ordinal at least 1 and at most $n-1$, so $\text{ord}_fR\in\{2,3,...,n\}$. Through above discussion we yield $\text{ord}_fP=1\Leftrightarrow P=(1)$.\\~\\
If $f$ is a simple standard rule with $|f|=n$, consider position $P=(n,1,1,1,...,1)$ where there are $(n-1)$ number of 1s. Such $\text{ord}_fP=1$ but obviously $P\neq(1)$. $\square$\\~\\
{\bfseries Corollary 3.13}\quad For non-standard rule $f$, $\overline{\text{Mov}}\{P|\text{ord}_fP=1\}=\mathcal{L}_1$.\\~\\
This corollary is obvious according to theorem 3.12, and it is really useful for the following section.
\section{Exchange Theorem}
\subsection{Theorem Proof}
Exchange theorem gives us another way to judge rule isomorphism. Different from the previous simple rule theorem, we look at two rules with same number of players.\\~\\
{\bfseries Theorem 4.1}\quad (Exchange Theorem) Let $n$-player rule $f=(\alpha_1\alpha_2\alpha_3...\alpha_n)$, and there exists $\alpha_i$ and $\alpha_j$ satisfying $|\alpha_i-\alpha_j|=1$. We could get a new rule $g$ if we swap $\alpha_i$ and $\alpha_j$. If $$\overline{\text{Mov}}\{P|\text{ord}_fP=i-1\}\bigcap\overline{\text{Mov}}\{P|\text{ord}_fP=j-1\}=\varnothing$$
then $f\sim g$. Here $\overline{\text{Mov}}\mathcal{S}=\bigcup\limits_{P\in\mathcal{S}}\overline{\text{Mov}}P$.\\~\\
\textit{proof.}\quad We need to prove $\text{ord}_fP=\text{ord}_gP$ for any $P\in\mathcal{L}$. Let position $R$ be the smallest ($|R|=k$ minimum) counterexample. Then $R$ has $k$ different subsequent moves $Q_1$, $Q_2$, ..., $Q_k$. For those subsequent moves, the ordinals are the same for either $f$ or $g$. For each $1\leqslant i\leqslant k$, player 1 will get a point of $\alpha_{\text{ord}_fQ_i+1}$ if $R\to Q_i$ is used. As $\text{ord}_fR\neq\text{ord}_gR$ and rule $f$ and $g$ are only different at two ordinals, there must exists two subsequent moves, say $Q_1$ and $Q_2$, such that $\alpha_{\text{ord}_fQ_1+1}=\alpha_i$, $\alpha_{\text{ord}_fQ_2+1}=\alpha_j$. Then, $\text{ord}_fQ_1=i-1$, $\text{ord}_fQ_2=j-1$. Note that $R\to Q_1$, $R\to Q_2$, so $$R\in\overline{\text{Mov}}\{P|\text{ord}_fP=i-1\}\bigcap\overline{\text{Mov}}\{P|\text{ord}_fP=j-1\}\neq\varnothing$$
which is a contradiction. So no such $R$ as counterexample exists. Therefore $f\sim g$.$\square$
\subsection{Applications}
With this new pivotal theorem in hand, we can explore new kinds of isomorphisms. \\~\\
{\bfseries Corollary 4.2}\quad For $n$-player rule $f=(\alpha_1\alpha_2...\alpha_n)$, if there exists $n\geqslant i>m_f>1$ such that $|\alpha_2-\alpha_i|=1$, then we can swap $\alpha_2$ and $\alpha_i$ and get an isomorph rule.\\~\\
\textit{proof.}\quad By corollary 3.10, corollary 3.12 and Exchange Theorem, as $$\overline{\text{Mov}}\{P|\text{ord}_fP=2-1\}\bigcap\overline{\text{Mov}}\{P|\text{ord}_fP=i-1\}\subset\mathcal{L}_1\bigcap\mathcal{L}_2=\varnothing$$
so swapping $\alpha_2$ and $\alpha_i$ can achieve isomorph rules.$\square$\\~\\
Furthermore, we can completely categorize all rules with no more than 4 players using only Simple Rule Theorem and Exchange Theorem:
\begin{itemize}
    \item All rules with the first value non-zero is non-simple by Simple Rule Theorem.
    \item 1-player rule: only one, (0).
    \item 2-player rule: (01), the standard 2-player chomp.
    \item 3-player rule: (012), (021).
    \item 4-player rule: (0123), $(0132)\sim(0231)^*$, $(0213)\sim(0312)^*$, $(0321)\sim(021)^{**}$.
    \item So, under equivalence of isomorphism, with no more than 4 players, there are only seven rules: (0), (01), (012), (021), (0123), (0132), (0213).
\end{itemize}
$*$: Exchange Theorem, or more easily, by corollary 4.2. \\
$**$: Consider position $P$ with $\text{ord}_{(0321)}P=4$, then $P\in\mathcal{L}_2$ as $m_{(0321)}=2<4$. Note that $P=(a_1, a_2, ..., a_k)\to (a_1)\in\mathcal{L}_1$ is always a legal move and so player 1 can get at least 2 points, so $\text{ord}_{(0321)}P=2$ or 3, which is a contradiction. Thus, for any position $P$, $\text{ord}_{(0321)}P\leqslant3<4$. By Simple Rule Theorem, $(0321)\sim(032)=(021)$. \\~\\
It is really interesting to investigate all the isomorphs for small rules.
\section{Further Questions}
\begin{enumerate}
    \item Algorithmic concerns: we know that judging the ordinal of a given position $P$ under rule (01) (a.k.a standard 2-player chomp) is NP-complete, but under other multiplayer rules, is it possible to solve for the ordinal of a random position within polynomial time?
    \item Rule classification: now that we have listed chomp rules with no more than 4 players. Can we categorize multiplayer chomp rules by isomorphism, by ascending number of players, as many as possible?
    \item If given rule $f\not\sim g$, we can find a position $P$ with least volume $V(f,g)$ such that $\text{ord}_fP\neq\text{ord}_gP$. Denote $V_n=\max\limits_{|f|,|g|\leqslant n}\{V(f,g)|f\not\sim g\}$, then what is $V_n$ with respect to $n$? $V_n$ will give us a general picture of how many positions are there to calculate if we want to classify all the rules.
\end{enumerate}
\section*{Acknowledgements}
I would like to thank Professor Doron Zeilberger and Young Investigator Hehui Wu for helpful comments that improved the manuscript.

\end{document}